\documentclass[]{interact}

\usepackage{epstopdf}% To incorporate .eps illustrations using PDFLaTeX, etc.
\usepackage{subfigure}% Support for small, `sub' figures and tables
\usepackage{float}

\usepackage[numbers,sort&compress]{natbib}% Citation support using natbib.sty
\bibpunct[, ]{[}{]}{,}{n}{,}{,}% Citation support using natbib.sty
\renewcommand\bibfont{\fontsize{10}{12}\selectfont}% Bibliography support using natbib.sty
\makeatletter% @ now becomes a letter
\def\NAT@def@citea{\def\@citea{\NAT@separator}}% Macro to suppress spaces between citations using natbib.sty
\makeatother% @ becomes a symbol again

\theoremstyle{plain}% Theorem-like structures provided by amsthm.sty
\newtheorem{theorem}{Theorem}[section]

\theoremstyle{definition}

\newtheorem{example}[theorem]{Example}

\theoremstyle{remark}

\usepackage[colorlinks=true,linkcolor=red, citecolor=blue, urlcolor=red]{hyperref}

\begin{document}

\articletype{}

\title{Asymmetric  Matrix-Valued Covariances for Multivariate Random Fields on Spheres}

\author{   \name{Alfredo Alegr\'ia\textsuperscript{1}\thanks{Alfredo Alegr\'ia. Email: alfredo.alegria.jimenez@gmail.com}, Emilio Porcu\textsuperscript{2,1}\thanks{Emilio Porcu. Email: georgepolya01@gmail.com} and Reinhard Furrer\textsuperscript{3}\thanks{Reinhard Furrer. Email: reinhard.furrer@math.uzh.ch}}  \affil{\textsuperscript{1}Department of Mathematics, University Federico Santa Maria,  Valparaiso, Chile;   \\ \textsuperscript{2}School of Mathematics and Statistics, University of Newcastle, Newcastle upon Tyne, United Kingdom;\\ \textsuperscript{3}Department of Mathematics and Department of Computer Science, University of Zurich, Zurich, Switzerland.}   }

\maketitle

\begin{abstract}
Matrix-valued covariance functions are crucial to geostatistical modeling of multivariate spatial data.  The classical assumption of symmetry of a multivariate covariance function is overlay restrictive and has been considered as unrealistic for most of real data applications. Despite of that, the literature on asymmetric covariance functions has been very sparse. In particular, there is some work related to asymmetric covariances on Euclidean spaces, depending on the Euclidean distance.  However, for data collected over large portions of planet Earth,  the most natural  spatial domain is a sphere,   with the corresponding geodesic distance being the natural metric.   In this work, we propose a strategy based on spatial rotations to  generate asymmetric covariances for multivariate random fields on the $d$-dimensional unit sphere.  We illustrate through simulations  as well as real data analysis that our proposal allows to achieve improvements in the predictive performance in comparison to the symmetric counterpart.
\end{abstract}

\begin{keywords}
 Cauchy model; Geodesic distance;  Global data; Rotation group; Wendland model.
\end{keywords}

\section{Introduction}

In recent years, the use of multivariate geostatistical models has become increasingly popular in  disciplines  such as  climatology, environmental sciences and oceanography.   Specifically, the    georeferenced observations  are assumed as a partial realization of a  multivariate  Gaussian random field (MGRF) and   the spatial dependencies between the observed data are  modeled through the  matrix-valued covariance function of the field.   Covariance functions are the main tool for  spatial interpolation methods and  the development of  accurate predictions. The  best linear unbiased predictor, under  a mean  squared error criterion,  is known in multivariate geostatistics as the co-kriging predictor.  We refer the reader to \cite{Wackernagel:2003} for a detailed study on MGRFs and \cite{Genton:Kleiber:2014} for a review  on parametric families of matrix-valued covariance functions.

A typical assumption in geostatistics is that of spatial isotropy, which means that the multivariate covariance depends solely on spatial distance.   The practical and theoretical contributions  in spatial statistics   have been mainly developed on  Euclidean spaces  with the corresponding Euclidean metric.   This choice is supported by the existence of several methods to project the geographical locations onto the plane.     However,    for  data recorded  over  large portions of planet Earth,  such projection generate significant distortions in the distances.  The problem of an appropriate metric has been discussed by \cite{BIOM:BIOM040302}, \cite{gneiting2013} and \cite{PBG16}.   Indeed, if we agree that a sphere is a good approximation of planet Earth,  it is more realistic  to work under the framework of random fields indexed on a sphere \citep{marinucci2011random}.  Moreover, the most natural metric to describe the spatial interactions is the geodesic distance, which corresponds to the shortest path joining two points on the spherical surface.

The construction of flexible families of matrix-valued covariance functions on spheres is a challenging problem, since the corresponding covariance matrices must be positive definite.  The literature in this direction has been very sparse with the work of   \cite{PBG16}, who derive  conditions for the validity of  the  multivariate Mat\'ern  and Cauchy models, being a notable exception. On the other hand, \cite{ma2016stochastic} shows that  compactly supported models, introduced on Euclidean spaces, can be  adopted on  specific low dimensional spheres.  Additional  approaches based on latent dimensions are explored by \cite{alegria2017covariance}. 

The parametric covariance families mentioned above depend solely on the geodesic distance (geodesic isotropy). In particular, spatial isotropy  implies symmetry of the matrix-valued covariance structure.    However,  in real data applications, asymmetrical behaviors are usually observed.  Indeed, we need to go beyond the class of geodesically isotropic covariances  to generate asymmetric models.  In this paper, we propose a strategy   based on spatial rotations to generate   asymmetric covariances from  symmetric models.  We show that the inclusion of asymmetry produce more flexibility  in the statistical modeling, and particularly, that the improvements in terms of predictive performance are significant. Our proposal is the spherical  counterpart of the work developed by \cite{ li:an} on Euclidean spaces.

The remainder of the article is as follows. Section \ref{background} contains the background material on MGRFs on spheres as well as a brief summary of some symmetric  covariance families, which are the building block to construct asymmetric models. Section \ref{sec_asy} introduces a general approach to produce asymmetric covariances.   In Section \ref{simulation}, we show through simulation experiments that asymmetric models can provide  better predictive results in comparison to the symmetric counterparts.   A bivariate data set is analyzed in Section \ref{data}. Finally, Section \ref{discussion} contains some conclusions.

\section{Background}
\label{background}

We briefly review the literature on MGRFs on spheres and the parametric families of matrix-valued covariance functions introduced in the literature. Although the models considered in this section are symmetric,  they are the building block to construct  explicit asymmetric models in the subsequent sections.

\subsection{MGRFs on Spheres}

Throughout the article, we consider  MGRFs, $\{\bm{Z}(\bm{x})=(Z_1(\bm{x}),\hdots,Z_p(\bm{x}))^\top: \bm{x}\in \mathbb{S}^d\}$,  defined on the $d$-dimensional unit sphere $\mathbb{S}^d=\{\bm{x}\in\mathbb{R}^{d+1} : \|\bm{x}\|=1 \}$, $d\in\mathbb{N}$,  where  $p\in\mathbb{N}$ denotes the number of components of the field, $^\top$  the transpose operator and $\|\cdot\|$ the Euclidean norm in $\mathbb{R}^{d+1}$.  We  call $\bm{Z}(\bm{x})$  a $p$-variate random field. The  natural  metric on the spherical surface is the geodesic distance defined through the mapping   $\theta:\mathbb{S}^d \times \mathbb{S}^d \rightarrow [0,\pi]$ given by $ \theta := \theta(\bm{x},\bm{y}) = \arccos(\bm{x}^\top\bm{y})$.

 We assume, without loss of generality,  that $\bm{Z}(\bm{x})$ has zero mean, since we are essentially interested on its covariance structure. 
 Let $i,j\in\{1,\hdots,p\}$  and  consider the matrix-valued mapping ${\bm{F}} : \mathbb{S}^d \times \mathbb{S}^d \rightarrow \mathbb{R}^{p\times p}$, with $(i,j)$-th element defined through
\begin{equation*}
  \label{eq_cross}
  {F}_{ij}(\bm{x},\bm{y}) := \text{cov}\{Z_i(\bm{x}),Z_j(\bm{y})\},   \qquad  \bm{x},\bm{y}\in\mathbb{S}^d.
  \end{equation*}
The diagonal elements, ${F}_{ii} :  \mathbb{S}^d \times \mathbb{S}^d \rightarrow \mathbb{R}$, are called  marginal covariances,  whereas the off-diagonal elements,  ${F}_{ij}: \mathbb{S}^d \times \mathbb{S}^d \rightarrow \mathbb{R}$, for $i\neq j$,  are called cross-covariances.   The mapping  ${\bm{F}}$   is (semi) positive definite, i.e.,  for all positive integer   $n$, $\{\bm{x}_1,\hdots,\bm{x}_n\} \subset \mathbb{S}^d $ and $\{\bm{a}_1,\hdots,\bm{a}_n\}\subset \mathbb{R}^p$, the following inequality holds
\begin{equation*} \label{posdef}
   \sum_{\ell=1}^{n} \sum_{r=1}^n   \bm{a}_{\ell}^\top { \bm{F}}(\bm{x}_\ell,\bm{x}_r)  \bm{a}_{r} \geq 0.
 \end{equation*}

 Note that, in general, the cross-covariance function  is not symmetric in the sense that $F_{ij}(\bm{x},\bm{y}) \neq F_{ji}(\bm{x},\bm{y})$ and $F_{ij}(\bm{x},\bm{y}) \neq F_{ij}(\bm{y},\bm{x})$. However,  when this covariance depends on $\bm{x}$ and $\bm{y}$ only through its geodesic distance (so-called geodesic isotropy), such a covariance function is clearly symmetric. The rest of this section focuses on geodesically isotropic models. 
 
 Following \cite{PBG16}, we call $\Psi_{d}^p$ the class of continuous mappings $\bm{C}: [0,\pi]\rightarrow \mathbb{R}^{p\times p}$, being continuous and such that, for a positive definite function ${\bm{F}} : \mathbb{S}^d \times \mathbb{S}^d \rightarrow \mathbb{R}^{p\times p}$, we have  ${\bm{F}}(\bm{x},\bm{y}) = \bm{C}(\theta(\bm{x},\bm{y}))$. The mapping  $\bm{C}$ is called radial part of a geodesically isotropic covariance ${\bm{F}}$. Moreover, the  inclusion $\Psi_{d+1}^p \subset \Psi_d^p$, for all positive integer $d$, is strict. 

Throughout, $\mathcal{P}_k^\lambda(\cdot)$ denotes the $\lambda$-ultraspherical polynomial of degree $k$ \citep{Abramowitz-Stegun:1965}, which is defined implicitly as
\begin{equation*}
 \frac{1}{(1+r^2-2r\mu)^\lambda}  =  \sum_{k=0}^\infty  r^k \mathcal{P}_k^\lambda(\mu), \qquad \mu\in[-1,1], r\in(-1,1).
\end{equation*}
  The following result \citep{hannan2009multiple, Yaglom:1987} characterizes completely the class  $\Psi_{d}^p$ in terms of ultra\-spherical polynomials. The equalities and summability conditions must be understood in a componentwise sense.

\begin{theorem} \label{teo}
 Let $d$ and $p$ be   positive integers. Let $\bm{C}:[0,\pi]\rightarrow \mathbb{R}^{p\times p}$ be a continuous mapping with elements $C_{ij}(\cdot)$, such that $C_{ii}(0)<+\infty$, for all $i=1,\hdots,p$. Then,  $\bm{C}$ is a member  of the class $\Psi_{d}^p$ if and only if it admits the  representation
\begin{equation*}
\label{sch_representation}
\bm{C}(\theta) = \sum_{k=0}^\infty {\bm{B}}_k \mathcal{P}_k^{(d-1)/2}(\cos\theta), \qquad \theta \in [0,\pi],
\end{equation*}
where $\{\bm{B}_k\}_{k=0}^\infty$ is a sequence of symmetric, positive definite and summable matrices of order $p\times p$.  
\end{theorem}

%We recall that the symmetry of the   matrices $\bm{B}_k$ in  Theorem \ref{teo} implies  that   $C_{ij}(\theta) = C_{ji}(\theta)$, for all $i,j=1,\hdots,p$.  Thus, in order to obtain asymmetric models, we need to evade from the class $\Psi_{d}^p$.

\subsection{Parametric Families in $\Psi_d^p$}
\label{review_cross_covariances}

We now provide a  review on  parametric families of geodesically isotropic  matrix-valued covariance functions. We focus on the Mat\'ern and Cauchy families, as well as on compactly supported models. In particular, the Mat\'ern and Cauchy families are valid  on spheres of any dimension, whereas the compactly supported models require additional restrictions on the dimension.

\begin{description}

\item[Mat\'ern model]  Let $K_\nu$ be the modified Bessel function of second kind and $\Gamma$ be the Gamma function \citep{Abramowitz-Stegun:1965}. The multivariate Mat\'ern model  on the sphere \citep{PBG16} is given by
\begin{equation*}
\label{matern}
C_{ij}(\theta) =  \sigma_i \sigma_j \rho_{ij} \frac{2^{1-\nu}}{\Gamma(\nu)}  \left(\frac{\theta}{c_{ij}} \right)^{\nu} K_{\nu}\left(\frac{\theta}{c_{ij}}\right), \qquad \theta\in[0,\pi],
\end{equation*} 
provided that  $\sigma_i$ and $c_{ij}$ are positive, $0<\nu\leq 1/2$,  $|\rho_{ij}|\leq 1$,  with $\rho_{ii} = 1$, for all $i,j=1,\hdots,p$.  Here, $\sigma_i$ is the standard deviation of the $i$-th component of the field, $c_{ij}$ are scale parameters, $\rho_{ij}$ is the collocated correlation coefficient between the components $i$ and $j$, and $\nu$ is a smoothness parameter.   Some additional conditions for the parameters are required (see \cite{PBG16}). In the following sections, we focus on  parsimonious bivariate ($p=2$)  models and such conditions are relaxed. For instance, the  parameterization $c_{12} = \max\{c_{11},c_{22}\}$ has associated  condition
\begin{equation}
\label{condition_matern}
|\rho_{12}| \leq   \left(\frac{c_{11}c_{22}}{c_{12}^2}\right)^{\nu}.
\end{equation} 
For a more detailed  study of the parametric conditions, we refer the reader to  \cite{gneiting2010matern},  \cite{apanasovich2012valid} and  \cite{PBG16}.
%In particular, if $\nu=1/2$, Equation (\ref{matern}) reduces to the exponential model
%\begin{equation}
%C_{ij}(\theta) = \sigma_i \sigma_j \rho_{ij}  \exp\left\{- \frac{3\theta}{ c_{ij}}\right\}, \qquad \theta\in[0,\pi].
%\end{equation} 

\item[Cauchy model]   The conditions for matrix-valued covariance functions of  generalized Cauchy type  on spheres have been developed by  \cite{PBG16}. The generalized Cauchy model is given by
\begin{equation*}
\label{cauchy}
C_{ij}(\theta) =  \sigma_i \sigma_j \rho_{ij}  \left(1 +  \frac{\theta^{\gamma}}{c_{ij}} \right)^{-\nu}, \qquad \theta\in[0,\pi],
\end{equation*} 
with   $\sigma_i$, $c_{ij}$  and $\nu$  being positive, $0<\gamma \leq 1$,   $|\rho_{ij}|\leq 1$ and $\rho_{ii} = 1$, for all $i,j=1,\hdots,p$.  Here, the interpretation of the parameters is similar to the Mat\'ern model.  Additional conditions must be imposed in order to obtain a positive definite mapping  (see \cite{PBG16}). For instance, in a bivariate case with $c_{12}=(c_{11}+c_{22})/2$, we have the restriction
\begin{equation}
\label{condition_cauchy}
\rho_{12}^2 \leq  \left( \frac{c_{11}c_{22}}{c_{12}^2} \right)^{\nu}.
\end{equation}
Moreover, \cite{PBG16} provide  explicit conditions on the parameters for the $p$-variate case, with $p> 2$.

\item[Compactly supported models]   \cite{ma2016stochastic} shows that a radial matrix-valued positive definite function on $\mathbb{R}^3$, with elements having compact support less than $\pi$, can be coupled with the geodesic distance  to obtain a geodesically isotropic positive definite function on $\mathbb{S}^3$.

 Formally, consider the mapping $\bm{\varphi}: [0,+\infty) \rightarrow \mathbb{R}^{p\times p}$, such that $(\bm{x},\bm{y}) \mapsto \bm{\varphi}(\| \bm{x}-\bm{y} \|)$, with elements $\varphi_{ij}(\cdot)$,  is positive definite on $\mathbb{R}^3$, with $\varphi_{ij}(r) = 0$, for all $r>\pi$ and $i,j=1,\hdots,p$. Therefore, the mapping $\bm{C}: [0,\pi] \rightarrow \mathbb{R}^{p\times p}$ defined by $\bm{C}(\theta) = \bm{\varphi}(\theta)|_{[0,\pi]}$ belongs to $\Psi_3^p$, where $f|_A$ denotes the restriction of $f$ to a Borel set $A\subset \mathbb{R}$.

In particular, the inclusion between the classes $\Psi_d^p$ imply that, under this construction principle, we can obtain models on spheres of dimension 1, 2 and~3. Examples can be generated from Table 1 in \cite{daley2015}. For instance, the Wendland model is defined as
 \begin{equation*}
\label{wg_sphere}
C_{ij}(\theta)= \sigma_i \sigma _j \rho_{ij} \left(1-\frac{\theta}{c_{ij}}\right)_{+}^{\nu} \left(1+\nu\frac{\theta}{c_{ij}}\right), \qquad \theta\in[0,\pi],
\end{equation*}
where $\nu\geq 2$, $0<c_{ij} \leq \pi$, $\sigma_i>0$, $|\rho_{ij}|\leq 1$ and $\rho_{ii} = 1$, for all $i,j=1,\hdots,p$. Here $(a)_{+} = \max\{a,0\}$ and the interpretation of the parameters is similar to the previous examples.   In addition, we have the condition (see \cite{daley2015})
\begin{equation}
\label{askey3}
\sum_{i\neq j} |\rho_{ij}|  (c_{ii}/c_{ij})^{\nu+1} \leq 1.
\end{equation}
In the subsequent sections, we consider a bivariate Wendland model with the parsimonious choice $c_{12}=\max\{c_{11},c_{22}\}$.

\end{description}

The characterization of the cross scale parameter $c_{ij}$ as a function of the marginal ones, $c_{ii}$ and $c_{jj}$, is a very useful strategy to avoid an excessive number of parameters in the estimation procedure. In some cases, this parsimonious strategy can provide better predictive results than the full model (see \cite{gneiting2010matern}).  

Additional parametric models can be constructed using latent dimensions  \citep{alegria2017covariance},  or by extending the  spherical convolution approach  \cite{ziegel2014convolution} as well as covariance approximation methods  \cite{jeong2015class} to the multivariate case.

\section{Asymmetric Covariances on $\mathbb{S}^d$}
\label{sec_asy}

A strategy  to construct asymmetric models is now provided. We start with a $p$-variate zero mean Gaussian field, $\bm{Z}(\bm{x})$, with geodesically isotropic, and thus  symmetric,  covariance structure  $\bm{F}(\bm{x},\bm{y}) = \bm{C}(\theta(\bm{x},\bm{y})) \in \Psi_d^p$.
  
Let $\{{\bm{R}}_{1}, \hdots, {\bm{R}}_p\} \subset \mathbb{R}^{(d+1)\times (d+1)}$ be a collection of rotation matrices, i.e., ${\bm{R}}_i$ is an orthogonal matrix with determinant one, for each $i=1,\hdots,p$.  Such matrices satisfy the relation ${\bm{R}}_i^\top = {\bm{R}}_i^{-1}$. Consider a new   $p$-variate random field  $\{\bm{Z}^a(\bm{x})=(Z_1^a(\bm{x}),\hdots,Z_p^a(\bm{x}))^\top: \bm{x}\in\mathbb{S}^d\}$ defined through
\begin{equation}
\label{rf2}
Z^a_i(\bm{x}) = Z_i({\bm{R}}_i \bm{x}),  \qquad  \bm{x}\in\mathbb{S}^{d},  \qquad i=1,\hdots,p.
\end{equation}
The field introduced in Equation (\ref{rf2})   is Gaussian and it has zero mean. Furthermore, the  covariance function associated to $\bm{Z}^a(\bm{x})$, denoted as  $\bm{F}^a : \mathbb{S}^{d} \times \mathbb{S}^{d}  \rightarrow \mathbb{R}^{p\times p}$,  is given by
\begin{eqnarray*}
\label{cov_asy}
{F}^a_{ij}(\bm{x},\bm{y})   &   :=   &   \text{cov}\{ Z^a_i(\bm{x}), Z^a_j(\bm{y}) \} \nonumber \\
 &    =    &   C_{ij}(  \theta(  {\bm{R}}_i \bm{x},  {\bm{R}}_j\bm{y} ))  \nonumber \\
 &  =  &   C_{ij}(  \theta(  \bm{x}, {\bm{R}}_i^{-1} {\bm{R}}_j\bm{y}  ) ), 
\end{eqnarray*}
for each $i,j = 1,\hdots,p$. Obviously, $\bm{Z}^a(\bm{x})$ preserves the marginal structure associated to $\bm{Z}(\bm{x})$,  since   ${F}^a_{ii}(\bm{x},\bm{y}) = C_{ii}(\theta(\bm{x},\bm{y}))$, for each $i=1,\hdots,p$.  However,  the cross-covariances  ${F}_{ij}^a$, for $i\neq j$,  are  not  necessarily   geodesically isotropic. Thus, $\bm{F}^a$ does not belong to the class $\Psi_d^p$.  Moreover,  we have that  ${F}_{ij}^a(\bm{x},\bm{y}) \neq {F}^a_{ji}(\bm{x},\bm{y})$, for $\bm{x},\bm{y}\in\mathbb{S}^d$.

Here, the rotation matrices are handled as additional parameters of the model. In order to avoid identifiability problems in the estimation of such parameters, we need to add some constraints in a similar fashion to the work of \cite{ li:an}. Since only the product ${\bm{R}}_i^{-1} {\bm{R}}_j$ is relevant, we propose the condition ${\bm{R}}_1 \cdots {\bm{R}}_p = {\bm{I}}_{d+1}$, where ${\bm{I}}_{d+1}$ denotes the identity matrix of order $(d+1)\times (d+1)$.  We illustrate these considerations for the cases $d=1$ and $d=2$.

\begin{example}

In the unit circle $\mathbb{S}^1$, we have the general rotation matrix $${\bm{R}}(\delta) : = \begin{bmatrix}
\cos\delta & \sin\delta\\
-\sin\delta & \cos\delta
\end{bmatrix}, \qquad   -\pi < \delta \leq \pi,$$  
where $\delta$ denotes the rotation angle. For each $i=1,\hdots,p$,  we define  ${\bm{R}}_i := {\bm{R}}(\delta_i)$.  
 Thus,  the constraint  ${\bm{R}}_1 \cdots {\bm{R}}_p ={\bm{R}}(\delta_1+\cdots+\delta_p)= {\bm{I}}_2$ is equivalent to the condition $(\delta_1+\cdots+\delta_p) \text{ mod } 2\pi = 0$.  In particular, we use the following identifiability condition 
\begin{equation}\label{IC}
\delta_1+\cdots+\delta_p = 0.
\end{equation}
Direct inspection shows that ${\bm{R}}_i^{-1} {\bm{R}}_j = {\bm{R}}(\delta_j-\delta_i)$  and $ \bm{x}^\top {\bm{R}}_i^{-1} {\bm{R}}_j\bm{y}  =  \cos(\theta)\cos(\delta_j-\delta_i)  + (x_1 y_2 - x_2 y_1) \sin(\delta_{j}-\delta_i)$, where $\theta := \theta(\bm{x},\bm{y})$ and $\bm{x}=(x_1,x_2)^\top,\bm{y}=(y_1,y_2)^\top \in\mathbb{S}^1$. 

   For instance, consider a bivariate scenario with $\delta_1=\eta/2$ and $\delta_2= - \eta/2$, for $-\pi < \eta<\pi$. In this case, we have the additional  asymmetry parameter $\eta$.  
   %Figure \ref{torus} depicts the covariance $(\bm{x},\bm{y}) \mapsto \exp\{  - \theta(\bm{x}, \bm{R}(\eta) \bm{y}) \}$, defined  on the torus $(\bm{x},\bm{y}) \in \mathbb{S}^1 \times \mathbb{S}^1$, with $\eta=0, 0.6, -0.6$.  Note that, by symmetry of the figure, we only need to plot  the covariance on the half of the toroidal surface. The  asymmetry parameter $\eta$ generates a shift effect in the covariance. This exponential structure is a particular example of the Mat\'ern model introduced in the previous sections, obtained when $\nu=1/2$.
    The general expression for the  bivariate asymmetric covariance functions on $\mathbb{S}^1$, generated from our approach, is  given by
$$   \bm{F}^a(\bm{x},\bm{y}) = \sum_{k=0}^\infty {\bm{B}}_k  \circ  \bm{P}_k(\bm{x},\bm{y}),$$
   where
   $$ \bm{P}_k(\bm{x},\bm{y})  =  \begin{pmatrix}   
\mathcal{P}_k^0(\cos\theta) & \mathcal{P}_k^0(\cos \theta \cos \eta  + (x_1 y_2 - x_2 y_1) \sin \eta )\\
\mathcal{P}_k^0(\cos \theta \cos \eta -   (x_1 y_2 - x_2 y_1) \sin \eta ) & \mathcal{P}_k^0(\cos\theta)
\end{pmatrix},$$
with $\circ$ denoting the Hadamard product. Recall that  $\{\bm{B}_k\}_{k=0}^\infty$ is a summable sequence of symmetric positive definite matrices. In this case, the ultraspherical polynomials reduces to the Tchebyshev polynomials $\mathcal{P}_k^0(\cos\theta) = \cos(k\theta)$.

%\begin{figure}
%\centering
%\caption{Exponential covariance function in dimension one, for different values of the asymmetry parameter $\eta$.}
%\label{torus}
%\includegraphics[scale=0.15]{asym1} \includegraphics[scale=0.15]{asym2} \includegraphics[scale=0.15]{asym3}
%\end{figure}

\end{example}

\begin{example}

Consider now the unit sphere $\mathbb{S}^2$. The \textit{Rodrigues's rotation formula} provides the following representation for a rotation matrix, with respect to an axis $\bm{\omega}=(\omega_1,\omega_2,\omega_3)^\top\in\mathbb{S}^2$ by an angle $\delta\in(-\pi,\pi]$,
\begin{equation*}
\label{rodrigues}
{\bm{R}}_{\bm{\omega}}(\delta) = {\bm{I}}_3 +(\sin\delta) \bm{\Omega} + (1-\cos\delta)\bm{\Omega}^2,  
\end{equation*} 
where  $$\bm{\Omega} = \begin{bmatrix} 0 & -\omega_3 & \omega_2\\
\omega_3 & 0 & -\omega_1\\
-\omega_2 & \omega_1 & 0 \end{bmatrix}.$$
Here, we take the collection of matrices ${\bm{R}}_{i} = {\bm{R}}_{\bm{\omega}}(\delta_i)$. Note that we have the relation ${\bm{R}}_i^{-1} {\bm{R}}_j = {\bm{R}}_{\bm{\omega}}(\delta_j-\delta_i)$.    For the bivariate case ($p=2$), considering the identifiability condition (\ref{IC}), $\delta_1=\eta/2$ and $\delta_2= -\eta/2$,   we have the additional parameters $(\eta,\bm{\omega}^\top)^\top \in (-\pi,\pi)\times\mathbb{S}^2$.  Since the pairs $(\eta,\bm{\omega}^\top)$ and $(-\eta,-\bm{\omega}^\top)$ define the same rotation, we must to restrict the range of $\eta$ to the interval $[0,\pi)$. The axis of rotation can be parameterized in terms of spherical coordinates: $\bm{\omega} = (\cos\alpha_1\sin\alpha_2, \sin\alpha_1\sin\alpha_2,\cos\alpha_2)^\top$, with $0 \leq \alpha_1 < 2\pi$ and $0\leq \alpha_2 \leq \pi$. In conclusion, we have the vector of additional parameters $(\eta,\alpha_1,\alpha_2)^\top$. 

\end{example}

\section{Simulation Study}
\label{simulation}

We illustrate  simulation experiments related to bivariate Gaussian  fields on $\mathbb{S}^2$ with asymmetric covariances. We use initial covariance functions either of Mat\'ern, Cauchy and Wendland type, as introduced in the previous sections.  In particular, we consider the mentioned models   with the following   parameterizations: 

\begin{itemize}

\item[{(M1)}] Mat\'ern model with $\nu=1/2$, which reduces to the exponential model
$$ C_{ij}(\theta) = \sigma_i\sigma_j \rho_{ij} \exp\left(   -  \frac{3\theta}{c_{ij}}  \right), \qquad \theta\in[0,\pi], \qquad  i,j=1,2,$$
where $c_{12}=\max\{c_{11},c_{22}\}$.
\item[{(M2)}] Cauchy model with $\nu=1$, 
$$ C_{ij}(\theta) = \sigma_i\sigma_j \rho_{ij}    \left(  1 +  \frac{19 \theta}{c_{ij}}  \right)^{-1},   \qquad \theta\in[0,\pi], \qquad i,j=1,2,$$
where $c_{12}=(c_{11}+c_{22})/2$.

\item[{(M3)}] Wendland model with $\nu = 4$,
$$ C_{ij}(\theta) = \sigma_i\sigma_j \rho_{ij}  \left( 1 - \frac{\theta}{c_{ij}}   \right)_{+}^4 \left(   1+ 4 \frac{\theta}{c_{ij}}  \right),   \qquad \theta\in[0,\pi],  \qquad i,j=1,2,$$
where $c_{12}=\max\{c_{11},c_{22}\}$.
\end{itemize}

\begin{figure}
\centering
\includegraphics[scale=0.25]{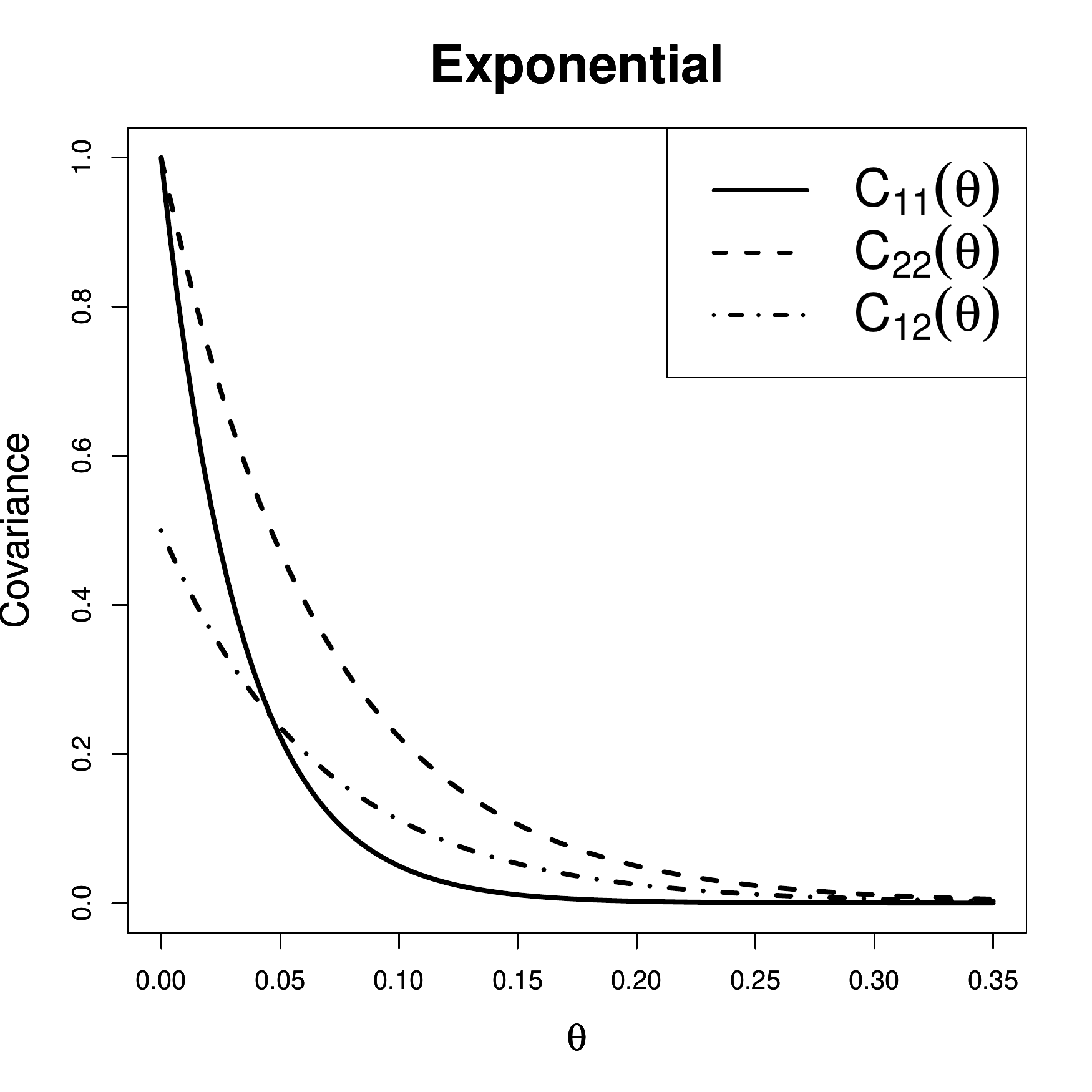} \includegraphics[scale=0.25]{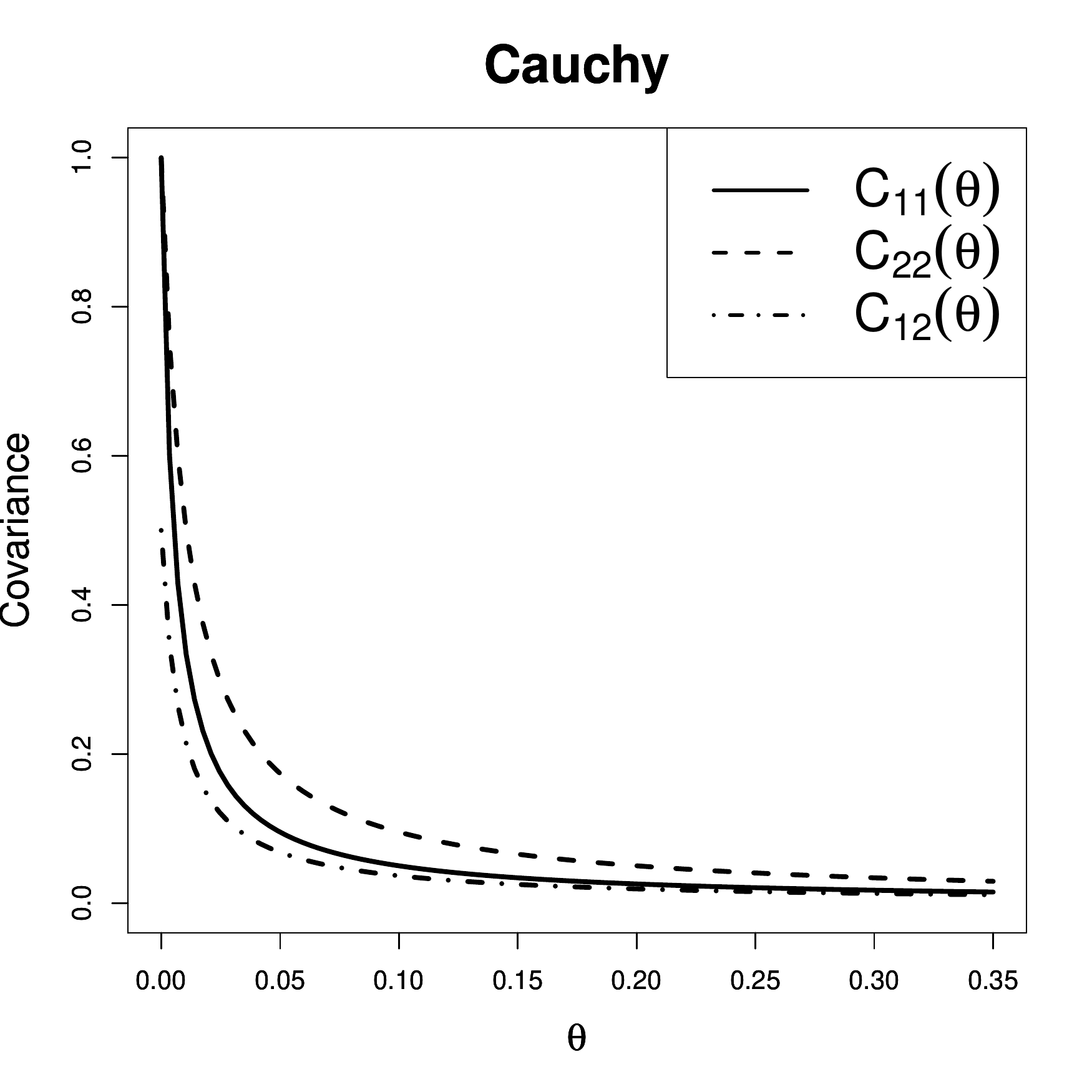} \includegraphics[scale=0.25]{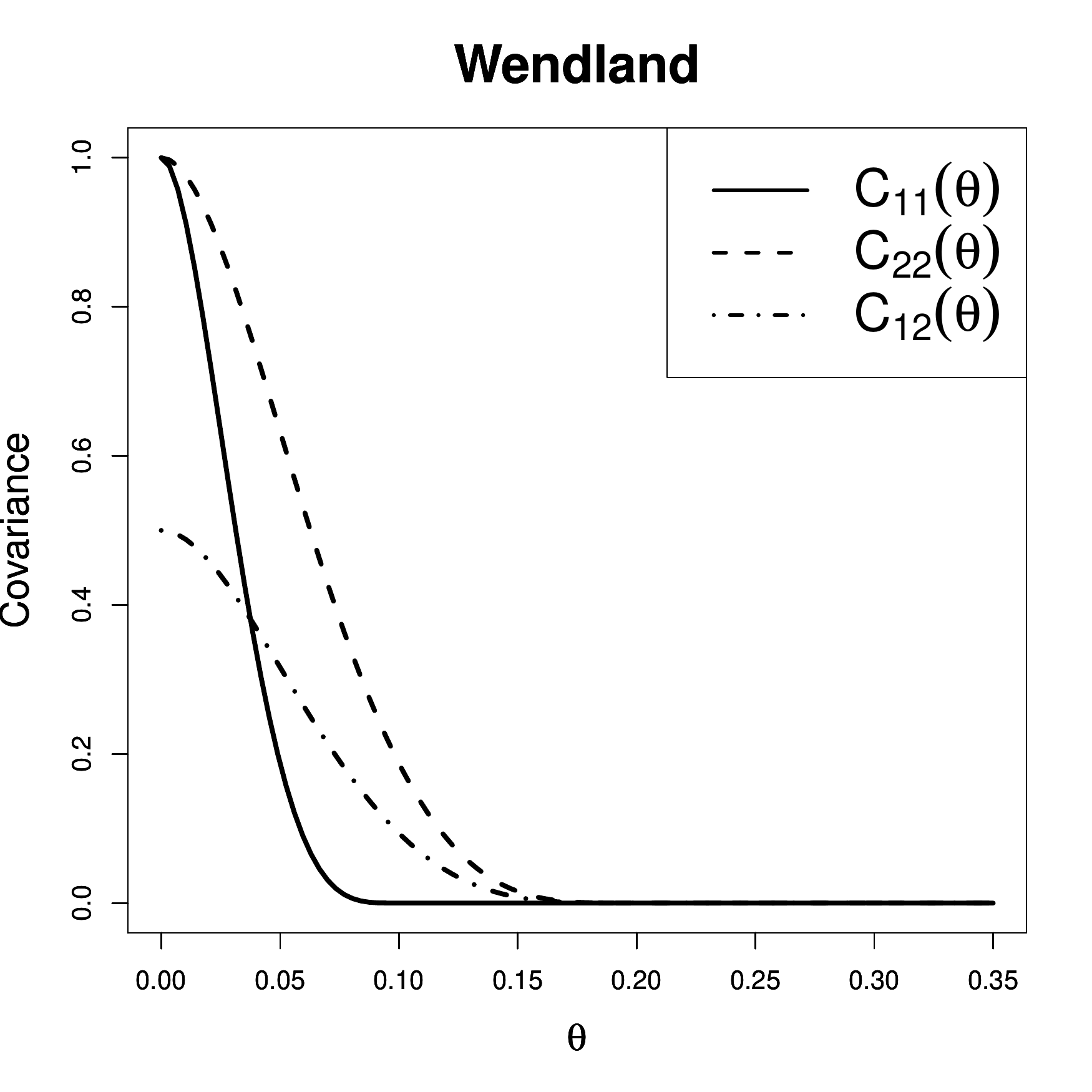}
\caption{Parametric families of matrix-valued covariances (M1)--(M3) considered in our simulation experiments. Here $\sigma_1^2 = \sigma_2^2 = 1$,  $\rho_{12}= 0.5$, $c_{11} = 0.1$, $c_{22} = 0.2$. }
\label{curves_cov}
\end{figure}

The parameterization of (M1) and (M2) ensures that $(C_{ij}(\theta)/C_{ij}(0))<0.05$ for $\theta>c_{ij}$. For model (M3),  $(C_{ij}(\theta)/C_{ij}(0))$  is exactly equal to 0 beyond the  cut-off distance $c_{ij}$.  The large number of parameters considered in our simulation experiments   supports the parsimonious parameterization used for $c_{12}$.  In particular, for the generalized Cauchy model  (M2)  we consider $c_{12}$ as the average of the marginal scale parameters, since under this assumption we have the condition (\ref{condition_cauchy}), which can be checked easily. The same happens for the other models (see conditions (\ref{condition_matern}) and (\ref{askey3})).   Figure \ref{curves_cov} shows these covariances under a specific parameter setting:  $\sigma_1^2 = \sigma_2^2 = 1$,  $\rho_{12}= 0.5$, $c_{11} = 0.1$, $c_{22} = 0.2$.

In the special case where $c_{ij}=c>0$, for all $i,j=1,2$, we say that the covariances are separable, otherwise they are non-separable.  We make simulation experiments for both cases. For each non-separable model, after adding asymmetry, the vector of parameters is given by $(\sigma_1^2,\sigma_2^2,\rho_{12},c_{11},c_{22},\eta,\alpha_1,\alpha_2)^\top$. In the separable case, on the other hand, the vector of parameters reduces to $(\sigma_1^2,\sigma_2^2,\rho_{12},c,\eta,\alpha_1,\alpha_2)^\top$.  In our studies, we consider 225 spatial sites in a grid generated combining 15 equispaced azimuthal  and polar angles. Specifically, the polar angles are taken as $2\pi (k-1)/15$, whereas the azimuthal angles as $\pi (k-1)/15$, for $k=1,\hdots,15$.

\begin{figure}
\centering
\includegraphics[scale=0.4]{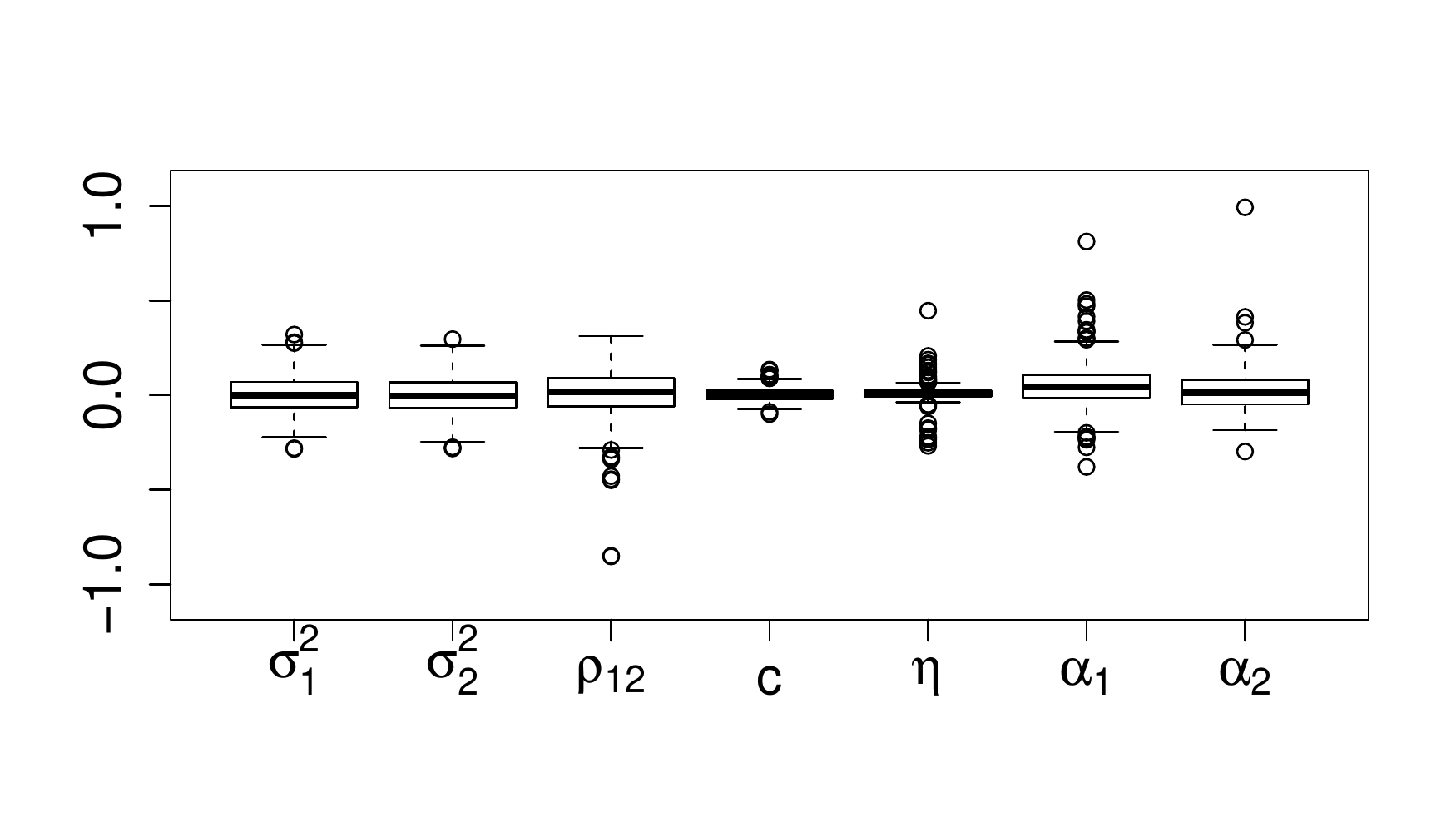}  \includegraphics[scale=0.4]{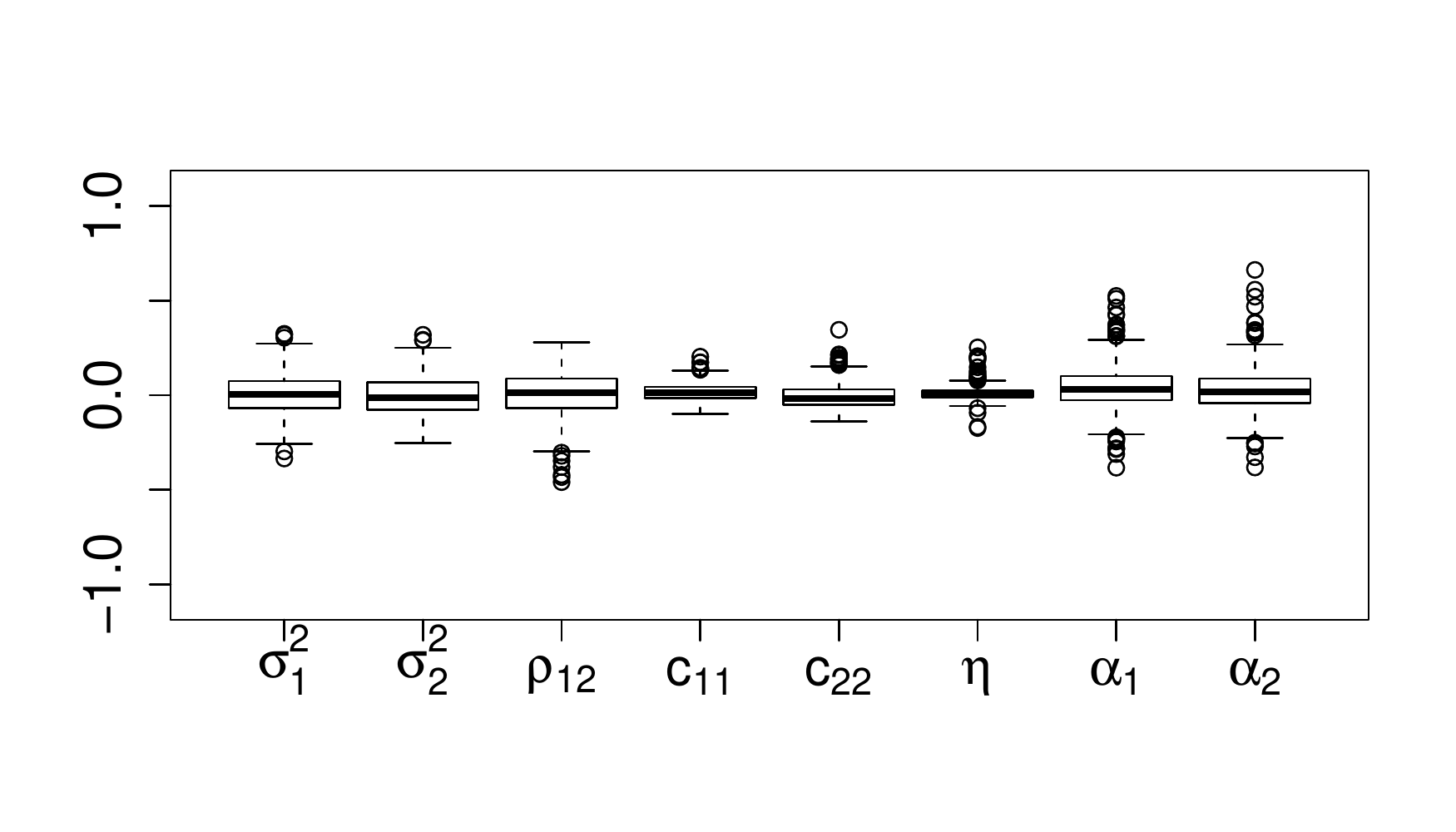}

\includegraphics[scale=0.4]{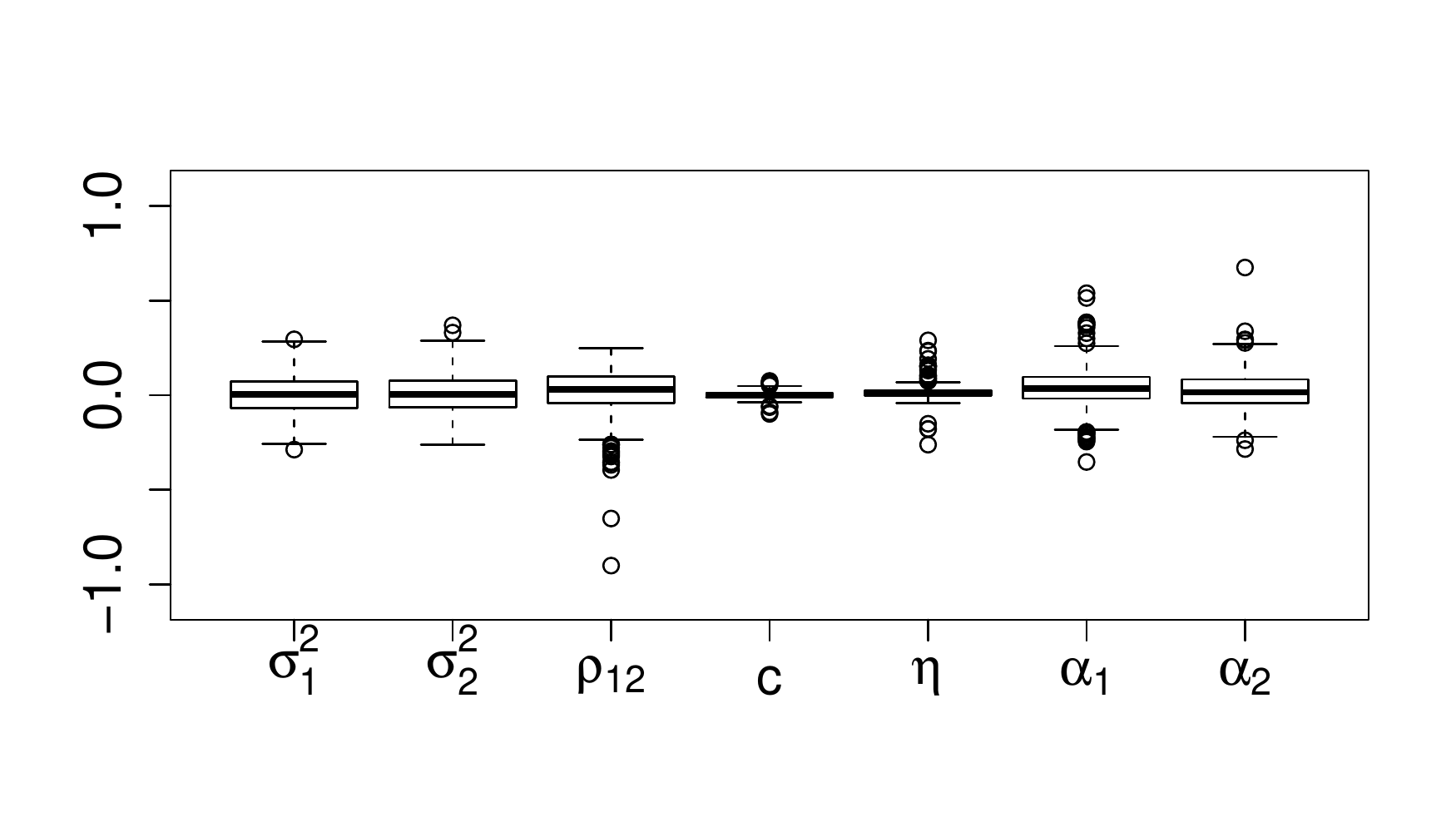} \includegraphics[scale=0.4]{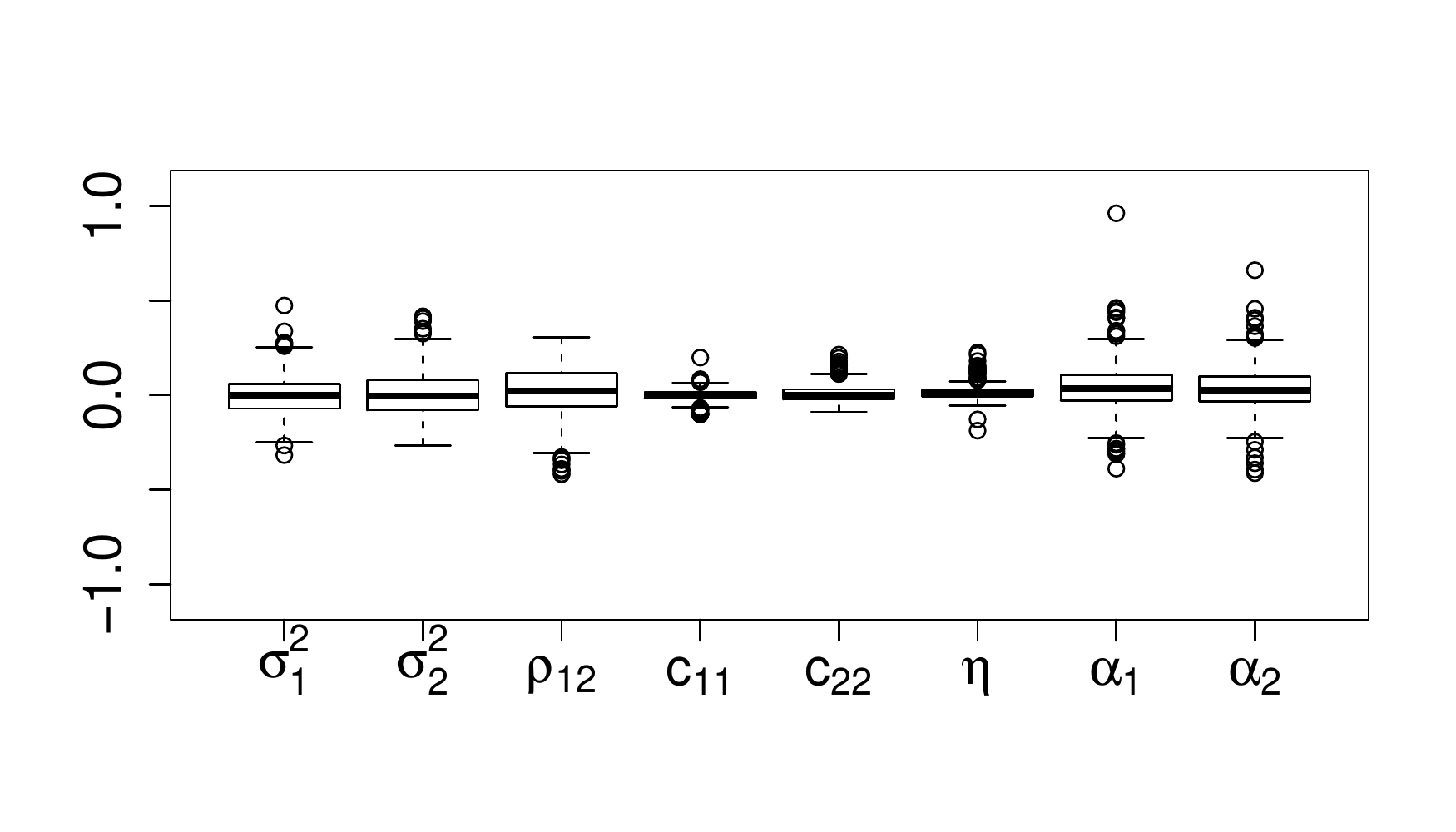}

\includegraphics[scale=0.4]{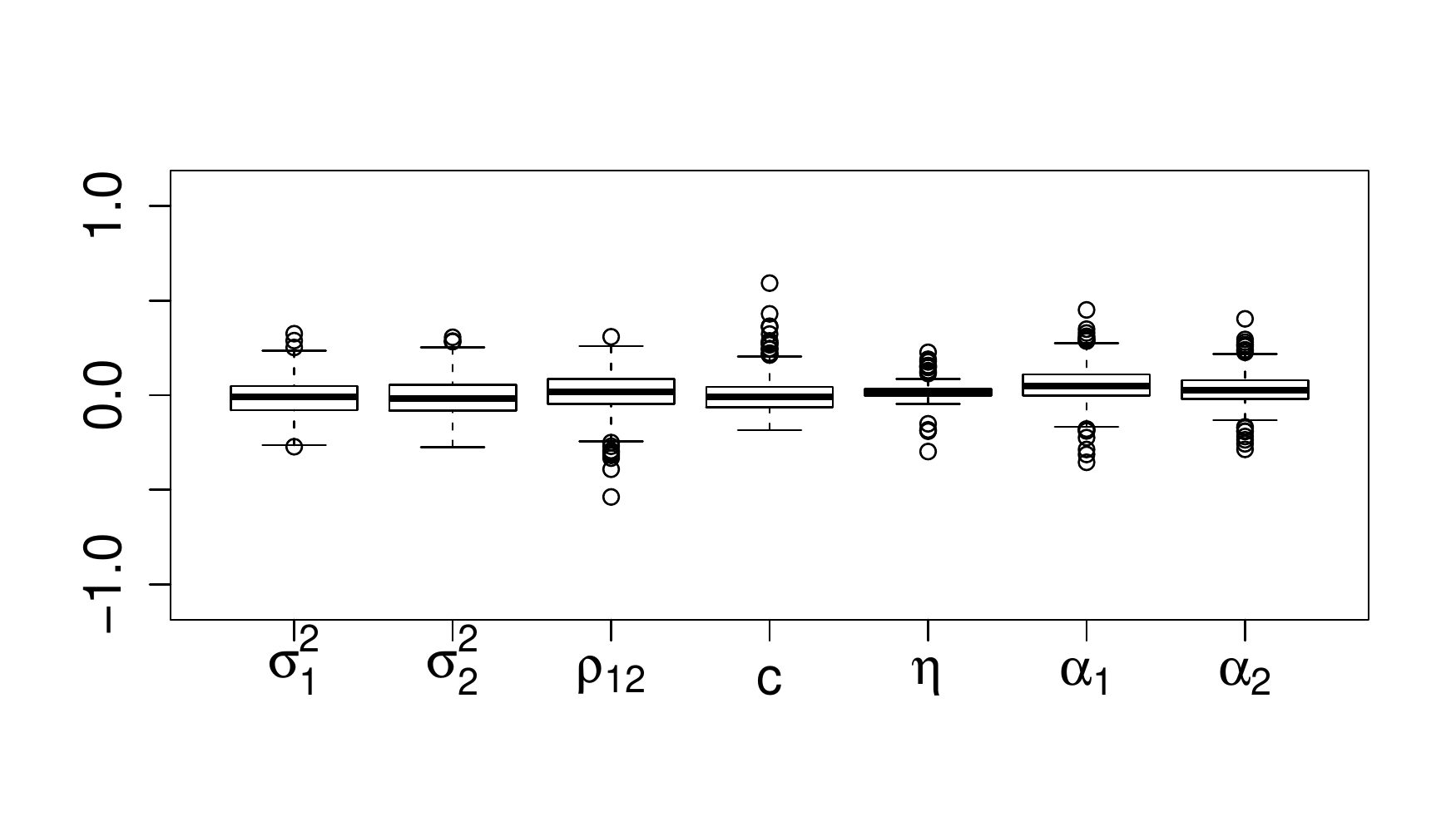} \includegraphics[scale=0.4]{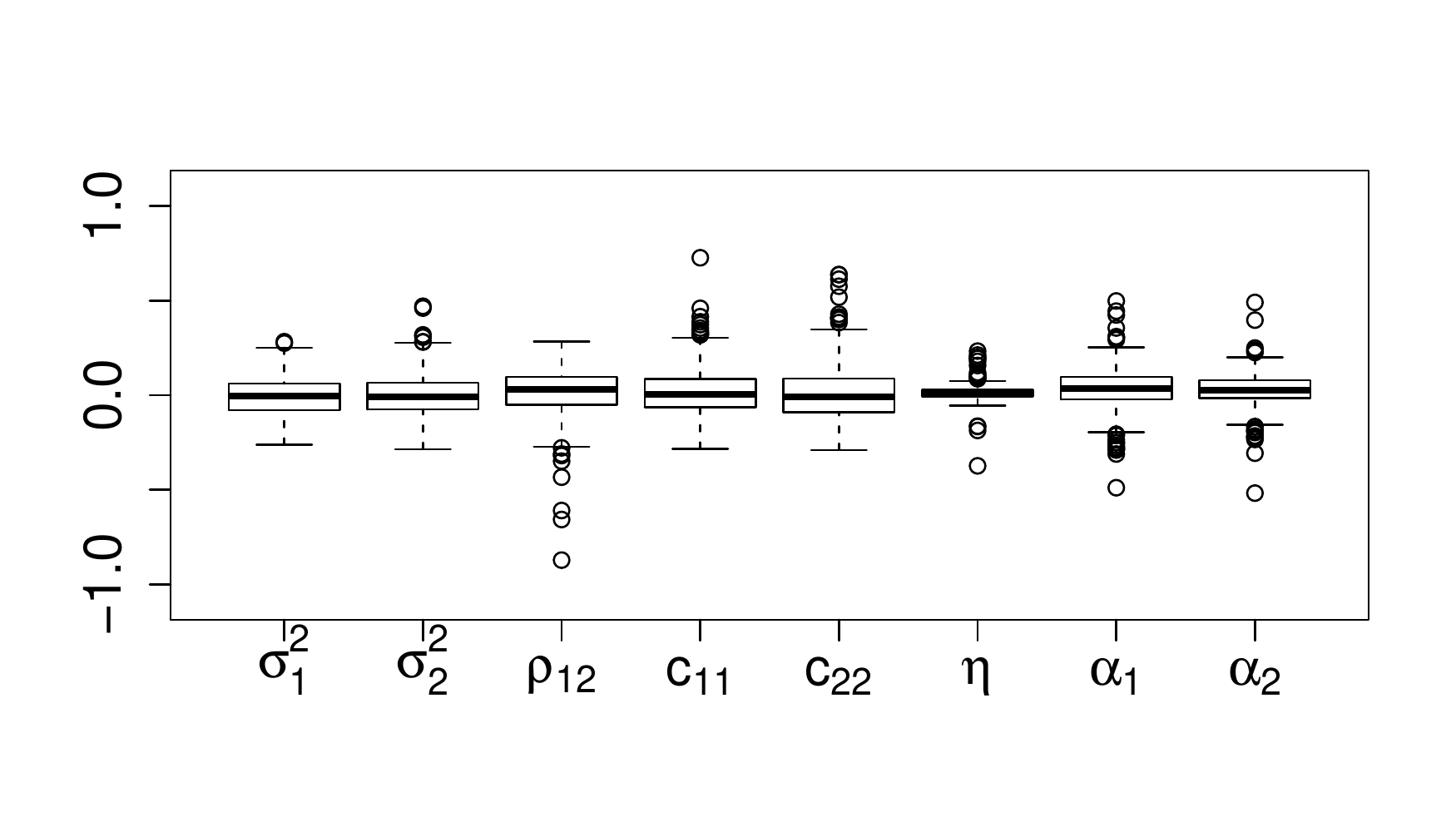} 
\caption{%Centered boxplots 
Bias of the pairwise CL estimates, on the basis of 500 independent repetitions. We consider models (M1)--(M3) (from  top to bottom), with  separable (left) and non-separable (right) structures. }
\label{boxplots}
\end{figure}

We assess the typical variance of the  likelihood-based estimates when the asymmetric effect is added. In particular, we consider the  pairwise composite likelihood (CL) approach developed by \cite{Bevilacqua2016}, which is a fast likelihood approximation method.   We set  $\sigma_1^2 = \sigma_2^2 = 1$,  $\rho_{12}= 0.5$, $c_{11} = 0.1$, $c_{22} = 0.2$, $\alpha_1=\alpha_2=\pi/2$ and $\eta=0.1$. Under the separable scenario, we set $c=0.1$. Figure \ref{boxplots} reports the bias %centered boxplots
 of the CL estimates for each model  on the basis of 500 independent simulated data sets. We have used a special case of  CL method, by considering  only pairs of data whose spatial distance is less than 1 radian, although in  our experience, smaller cut-off distances can also ensure stability of the method.  Even though we are dealing with a large number of parameters, the boxplots show that the direct maximization of the CL objective function provides quite reasonable results.  The empirical  $95\%$ confidence interval for $\eta$ is $[0.055, 0.190]$ in the non-separable model (M1), and similar intervals are obtained for the other models.   Another appealing estimation procedure   is to consider a two steps maximization routine (see \cite{ li:an}), which is useful in terms of computational complexity.

\begin{table}
 \caption{Cross-validation scores for models (M1)--(M3), separable and non-separable, using  symmetric (S) and asymmetric (A) versions, with $\rho_{12}=0.25$.}
 \label{mspe_simulation0}
\centering
\begin{tabular}{cccccccccccc} \hline 
\multicolumn{3}{c}{ }               &   \multicolumn{4}{c}{Separable}      &         &  \multicolumn{4}{c}{Non-separable}    \\ 
\cline{4-7} \cline{9-12}
\multicolumn{3}{c}{ }               &   \multicolumn{2}{c}{$\eta=0.1$}      & \multicolumn{2}{c}{$\eta=0.6$}  &       &  \multicolumn{2}{c}{$\eta=0.1$}  & \multicolumn{2}{c}{$\eta=0.6$}   \\ \cline{4-7} \cline{9-12}
              & &                      &        S              &         A                    &         S         &      A              &                  &  S      &  A           &         S         &      A                              \\  \hline
      (M1)    &     & MSPE      &    0.983          &   0.968                 &  0.964      & 0.935                    &                &  0.964  &    0.953         &   0.954       & 0.924                         \\
              &&  LSCORE         &     1.410         &  1.402                  &   1.398     & 1.383                   &                &  1.396  &    1.386         &   1.387       & 1.373                     \\ \hline
      %        && CRPS              &     2.112          &  2.096                  &  2.081      &  2.060                  &                &  2.068  &    2.058         &   2.071       & 2.043                       \\       \hline

      (M2)     & & MSPE          &     0.969         &  0.956                 &  0.961     & 0.930                     &              &  0.924  &    0.910          &   0.929       & 0.893                  \\
             && LSCORE            &     1.398         &  1.393                 &  1.395     & 1.381                  &                 &  1.357  &    1.348          &   1.352       & 1.327                       \\\hline
   %          && CRPS                 &     2.078         &  2.062                 &  2.083     &  2.057                 &                 &  1.997  &    1.989          &   2.026       & 1.9934                       \\       \hline
             
     (M3)     & & MSPE      &    0.981            &  0.966                     &  0.980        & 0.959                  &                &  0.993  &    0.988          &   0.972       & 0.945                 \\
             && LSCORE       &     1.409           &  1.399                     &  1.407        & 1.396                 &                 & 1.414   &    1.409          &   1.405       & 1.388                  \\\hline
%             && CRPS            &     2.099           &  2.080                     &  2.103        &  2.079                &                 &  2.117  &    2.113          &   2.093       & 2.075                    \\       \hline 
 \end{tabular}
 \end{table}

\begin{table}
 \caption{Cross-validation scores for models (M1)--(M3), separable and non-separable, using  symmetric (S) and asymmetric (A) versions, with $\rho_{12}=0.5$.}
 \label{mspe_simulation}
\centering
\begin{tabular}{cccccccccccc} \hline 
\multicolumn{3}{c}{ }               &   \multicolumn{4}{c}{Separable}      &         &  \multicolumn{4}{c}{Non-separable}    \\ 
\cline{4-7} \cline{9-12}
\multicolumn{3}{c}{ }               &   \multicolumn{2}{c}{$\eta=0.1$}      & \multicolumn{2}{c}{$\eta=0.6$}  &       &  \multicolumn{2}{c}{$\eta=0.1$}  & \multicolumn{2}{c}{$\eta=0.6$}   \\ \cline{4-7} \cline{9-12}
              & &                   &        S              &         A                  &         S         &      A              &                  &  S      &  A           &         S         &      A                              \\  \hline
      (M1)    &     & MSPE      &    0.969        &   0.945                &  0.979      & 0.919                    &            &  0.940  &    0.919         &   0.950       & 0.882                         \\
              &&  LSCORE     &     1.401         &  1.389                &   1.405       & 1.369                &                    &  1.383  &    1.372          &   1.387       & 1.345                     \\  \hline
    %          && CRPS          &     2.079         &  2.052                 &  2.096      &  2.017                  &                    &  2.050  &    2.021          &   2.062       & 1.987                       \\       \hline

      (M2)     & & MSPE          &     0.971        &  0.950                 &  0.975     & 0.908                     &             &  0.938  &    0.904          &   0.924       & 0.856                  \\
             && LSCORE        &     1.400           &  1.390                  &  1.406      & 1.366                 &                      &  1.356  &    1.339            &   1.358       & 1.326                       \\\hline
       %      && CRPS            &     2.083            &  2.052                 &  2.093       &  2.008               &                       &  2.043  &    1.995              &   2.029       & 1.949                       \\       \hline
             
     (M3)     & & MSPE      &    0.992          &  0.969                  &  0.993        & 0.911                      &             &  0.981  &    0.951          &   0.970       & 0.897                 \\
             && LSCORE  &     1.415           &  1.405                   &  1.414           & 1.364                 &                      & 1.408  &   1.393            &   1.403       & 1.366                  \\\hline
          %   && CRPS       &     2.111           &  2.078                   &  2.111           &  2.018                &                    &  2.087  &    2.064           &   2.100       & 2.006                    \\       \hline 
 \end{tabular}
 \end{table}

Now, we study the   predictive performance of  asymmetric models through a cross-validation analysis.   We simulate  bivariate fields   from the asymmetric versions of  models (M1)--(M3), with both separable and non-separable parameterizations. We set  $\sigma_1^2 = \sigma_2^2 = 1$, $c_{11}=0.1$, $c_{22}=0.2$, $\alpha_1=\alpha_2=\pi/2$, as in the previous example, and consider two scenarios for the collocated correlation coefficient: $\rho_{12}=0.25, 0.5$, and  two  scenarios for the asymmetry parameter: $\eta =   0.1, 0.6$. Again, in the separable case, we  consider a similar experiment with $c=0.1$.   We quantify the accuracy of the co-kriging predictor, with both the symmetric and asymmetric covariances,  in terms of the mean squared prediction error (MSPE) and  the log-score (LSCORE)  (see \cite{zhang2010kriging}). These indicators are evaluated using a drop-one prediction strategy. Tables \ref{mspe_simulation0} and  \ref{mspe_simulation} report the results for $\rho_{12}=0.25$ and $\rho_{12}=0.5$, respectively. It is clearly visible that the asymmetric models produce better results. Note that if  the asymmetry and the correlation between the fields are strong, then the improvements are more significant.

\section{A Real Data Example}
\label{data}

This section assesses the statistical performance of asymmetric models on a real data example. We analyze a bivariate data set of Temperatures (T) and Precipitations (P). These data outputs have been obtained from the  Community Climate System Model (CCSM4.0, \cite{doi:10.1175/2011JCLI4083.1}) provided by National Center for Atmospheric Research (NCAR), CO, USA.  The joint modeling of these variables is of major interest in climatological disciplines.

\begin{figure}
\centering
\includegraphics[scale=0.25]{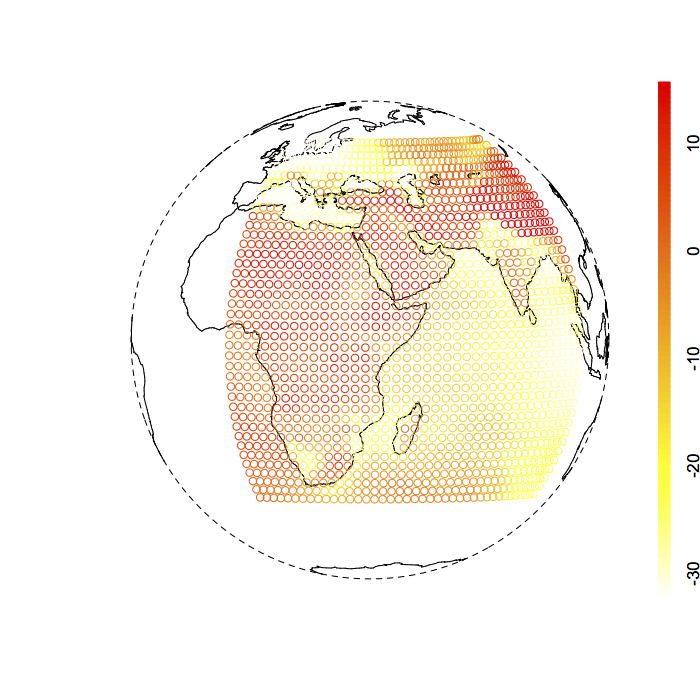} \includegraphics[scale=0.25]{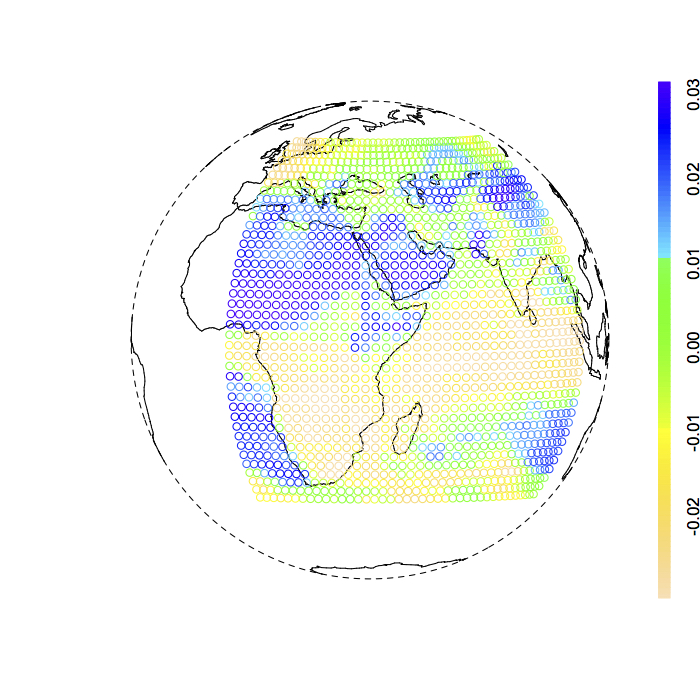}
\caption{Residuals for temperature (left) and cubic root of precipitation (right) of model output of December  2015.}
\label{datos_biv}
\end{figure}

The spatial resolution is 2.5 degrees in latitude and longitude.   We assume planet Earth as a sphere of radius 6378 kilometers. We pay attention to the geographical region delimited by the longitudes $[0,100]$ and latitudes $[-40,60]$ degrees (see Figure \ref{datos_biv}) and consider December 2015. The maximum geodesic distance between the points in this region is 70\% of the maximum distance on planet Earth. %Indeed, we are working over a global scale (see Figure \ref{datos_biv}).
We work with the cubic root of precipitations in order to attenuate the skewness of this variable.  The units are Kelvin degrees for temperatures and centimeters for transformed precipitations. In order to achieve geodesic isotropy, we remove through splines the trend and cyclic patterns of the variables, using  longitude and latitude as covariates. The resulting residuals are approximately Gaussian with zero mean.

We consider  the bivariate exponential model (M1) with the following features:
\begin{itemize}
\item Model 1. Separable exponential model without  asymmetry.
\item Model 2.  Non-separable exponential model  without  asymmetry.
\item Model 3. Separable exponential model after  adding  asymmetry.
\item Model 4.  Non-separable exponential model after adding  asymmetry.
\end{itemize}

Table \ref{estimaciones_datos} illustrates the pairwise CL estimates for each model.   We  have  considered  CL method with pairs of data whose geodesic distance is less than 1000 kilometers. The estimates of  the scale parameters are given in kilometers. We compare  the models in terms of the Log-Composite Likelihood (Log-CL) value at the optimum, as well as through their predictive performance.  The accuracy of the prediction, for each model, is evaluated in terms of MSPE and LSCORE,  using a drop-one prediction strategy (see Table \ref{scores_datos}).  It is clear that the asymmetric versions provide better results with respect to the symmetric ones.  Note that for this specific data set, Model 3 outperforms Model 2, which means that the inclusion of asymmetry can produce more flexibility in comparison to the inclusion of non-separability.

\begin{table}
\centering
\caption{Composite Likelihood estimates for Models 1-4.}
\label{estimaciones_datos}
\begin{tabular}{lccccccccc}  \hline  \raisebox{16pt}{~}
& $\widehat{\sigma}_{TT}^2$   &    $\widehat{\sigma}_{PP}^2$   &   $\widehat{\rho}_{TP}$   &   $\widehat{c}_{TT}$  & $\widehat{c}_{PP}$   &   $\widehat{c}_{TP}$   &  $\widehat{\eta}$   &  $\widehat{\alpha}_1$  & $\widehat{\alpha}_2$  \\ \hline
Model 1  \raisebox{16pt}{~}& 34.643  &  $9.262\times10^{-5}$   &  0.342  &   5102   &  5102   &  5102   &    -  &  -  &  - \\ 
Model 2 & 34.275  &  $9.363\times10^{-5}$   &  0.342  &   4789   & 5414   &  5414    &    -  &  -   &  - \\ 
Model 3 & 34.524  &  $9.235\times10^{-5}$    &  0.342  &  5089   & 5089    & 5089    &   0.043 &   1.512  &   2.678  \\
Model 4 & 34.098  &  $9.439\times10^{-5}$    &  0.342   &  4815  &  5427   & 5427    &  0.047  &  1.504    &  2.236 \\
  \hline 
\end{tabular}
\end{table}

\begin{table}
\centering
\caption{Comparison of Log-CL values and  cross-validation scores for Models 1-4.}
\label{scores_datos}
\begin{tabular}{cccccccccc}  \hline  
               &  $\#$ parameters &  Log-CL   &  MSPE   &  LSCORE     \\ \hline
Model 1   & 4 &    20460    &  1.309    &  $-1.235$           \\
Model 2   & 5 &    20485     &  1.306     &  $-1.236$             \\
Model 3   & 7 &    20500     &  1.157     &  $-1.271$             \\
Model 4   & 8 &    20518     &  1.139     &  $-1.277$           \\
  \hline 
\end{tabular}
\end{table}

\section{Conclusions}
\label{discussion}

We have provided a strategy to construct asymmetric matrix-valued covariance functions on spheres. We show that our proposal can produce significant improvements in the predictive performance of the model. Our developments have been illustrated through simulated data, under several models and parametric settings. We have also analyzed a real data set of temperatures and precipitations on a large portion of planet Earth, where the asymmetric models produce better results.

An interesting research direction is the construction of  parametric tests, based on  the asymmetry parameter $\eta$, to contrast the hypothesis of symmetry in a given data set. In particular, a challenging problem is the adaptation of the methodology used by  \cite{li2008testing}, where the authors  evaluate  several types of common assumptions on multivariate covariance functions, to the spherical context.

\section*{Acknowledgments}

Alfredo Alegr\'ia is supported by Beca CONICYT-PCHA/Doctorado Nacional/2016-21160371. Emilio Porcu is supported by Proyecto Fondecyt Regular number 1130647. Reinhard Furrer acknowledges support of the Swiss National Science Foundation SNSF-144973.  We additionally acknowledge the World Climate Research Programme's Working Group on Coupled Modelling, which is responsible for Coupled Model Intercomparison
Project (CMIP).

%\section*{Acknowledgments}
%The authors are very grateful to .

%\appendix
%\section*{Appendix}
%\addcontentsline{toc}{section}{Appendices}
%\renewcommand{\thesubsection}{\Alph{subsection}}
%\renewcommand{\theequation}{\Alph{subsection}.\arabic{equation}}

%\bibliographystyle{tfnlm}
\bibliographystyle{apalike}
\bibliography{mybib}

\begin{thebibliography}{}

\bibitem[Abramowitz and Stegun, 1965]{Abramowitz-Stegun:1965}
Abramowitz, M. and Stegun, I.~A. (1965).
\newblock Handbook of mathematical function: with formulas, graphs and
  mathematical tables.
\newblock Dover Publications.

\bibitem[Alegr{\'\i}a et~al., 2017]{alegria2017covariance}
Alegr{\'\i}a, A., Porcu, E., Furrer, R., and Mateu, J. (2017).
\newblock Covariance {F}unctions for {M}ultivariate {G}aussian {F}ields
  {E}volving {T}emporally {O}ver {P}lanet {E}arth.
\newblock {\em arXiv preprint}, arXiv:1701.06010.

\bibitem[Apanasovich et~al., 2012]{apanasovich2012valid}
Apanasovich, T.~V., Genton, M.~G., and Sun, Y. (2012).
\newblock A valid {M}at{\'e}rn class of cross-covariance functions for
  multivariate random fields with any number of components.
\newblock {\em Journal of the American Statistical Association},
  107(497):180--193.

\bibitem[Banerjee, 2005]{BIOM:BIOM040302}
Banerjee, S. (2005).
\newblock On geodetic distance computations in spatial modeling.
\newblock {\em Biometrics}, 61(2):617--625.

\bibitem[Bevilacqua et~al., 2016]{Bevilacqua2016}
Bevilacqua, M., Alegria, A., Velandia, D., and Porcu, E. (2016).
\newblock Composite {L}ikelihood {I}nference for {M}ultivariate {G}aussian
  {R}andom {F}ields.
\newblock {\em Journal of Agricultural, Biological, and Environmental
  Statistics}, 21(3):448--469.

\bibitem[Daley et~al., 2015]{daley2015}
Daley, D., Porcu, E., and Bevilacqua, M. (2015).
\newblock Classes of compactly supported covariance functions for multivariate
  random fields.
\newblock {\em Stochastic Environmental Research and Risk Assessment},
  29(4):1249--1263.

\bibitem[Gent et~al., 2011]{doi:10.1175/2011JCLI4083.1}
Gent, P.~R., Danabasoglu, G., Donner, L.~J., Holland, M.~M., Hunke, E.~C.,
  Jayne, S.~R., Lawrence, D.~M., Neale, R.~B., Rasch, P.~J., Vertenstein, M.,
  Worley, P.~H., Yang, Z.-L., and Zhang, M. (2011).
\newblock The community climate system model version 4.
\newblock {\em Journal of Climate}, 24(19):4973--4991.

\bibitem[Genton and Kleiber, 2015]{Genton:Kleiber:2014}
Genton, M.~G. and Kleiber, W. (2015).
\newblock Cross-{C}ovariance {F}unctions for {M}ultivariate {G}eostatistics.
\newblock {\em Statistical Science}, 30(2):147--163.

\bibitem[Gneiting, 2013]{gneiting2013}
Gneiting, T. (2013).
\newblock Strictly and non-strictly positive definite functions on spheres.
\newblock {\em Bernoulli}, 19(4):1327--1349.

\bibitem[Gneiting et~al., 2010]{gneiting2010matern}
Gneiting, T., Kleiber, W., and Schlather, M. (2010).
\newblock Mat{\'e}rn cross-covariance functions for multivariate random fields.
\newblock {\em Journal of the American Statistical Association},
  105(491):1167--1177.

\bibitem[Hannan, 2009]{hannan2009multiple}
Hannan, E. (2009).
\newblock {\em Multiple Time Series}.
\newblock Wiley Series in Probability and Statistics. Wiley.

\bibitem[Jeong and Jun, 2015]{jeong2015class}
Jeong, J. and Jun, M. (2015).
\newblock A class of mat{\'e}rn-like covariance functions for smooth processes
  on a sphere.
\newblock {\em Spatial Statistics}, 11:1--18.

\bibitem[Li et~al., 2008]{li2008testing}
Li, B., Genton, M.~G., and Sherman, M. (2008).
\newblock Testing the covariance structure of multivariate random fields.
\newblock {\em Biometrika}, pages 813--829.

\bibitem[Li and Zhang, 2011]{li:an}
Li, B. and Zhang, H. (2011).
\newblock An approach to modeling asymmetric multivariate spatial covariance
  structures.
\newblock {\em Journal of Multivariate Analysis}, 102(10):1445--1453.

\bibitem[Ma, 2016]{ma2016stochastic}
Ma, C. (2016).
\newblock Stochastic representations of isotropic vector random fields on
  spheres.
\newblock {\em Stochastic Analysis and Applications}, 34(3):389--403.

\bibitem[Marinucci and Peccati, 2011]{marinucci2011random}
Marinucci, D. and Peccati, G. (2011).
\newblock {\em Random fields on the sphere: representation, limit theorems and
  cosmological applications}, volume 389.
\newblock Cambridge University Press.

\bibitem[Porcu et~al., 2016]{PBG16}
Porcu, E., Bevilacqua, M., and Genton, M.~G. (2016).
\newblock Spatio-{T}emporal {C}ovariance and {C}ross-{C}ovariance {F}unctions
  of the {G}reat {C}ircle {D}istance on a {S}phere.
\newblock {\em Journal of the American Statistical Association},
  111(514):888--898.

\bibitem[Wackernagel, 2003]{Wackernagel:2003}
Wackernagel, H. (2003).
\newblock {\em Multivariate {G}eostatistics: An Introduction with
  Applications}.
\newblock Springer, New York, 3rd edition.

\bibitem[Yaglom, 1987]{Yaglom:1987}
Yaglom, A.~M. (1987).
\newblock {\em Correlation {T}heory of {S}tationary and {R}elated {R}andom
  {F}unctions. Volume {I}: {B}asic {R}esults}.
\newblock Springer, New York.

\bibitem[Zhang and Wang, 2010]{zhang2010kriging}
Zhang, H. and Wang, Y. (2010).
\newblock Kriging and cross-validation for massive spatial data.
\newblock {\em Environmetrics}, 21(3-4):290--304.

\bibitem[Ziegel, 2014]{ziegel2014convolution}
Ziegel, J. (2014).
\newblock Convolution roots and differentiability of isotropic positive
  definite functions on spheres.
\newblock {\em Proceedings of the American Mathematical Society},
  142(6):2063--2077.

\end{thebibliography}

\end{document}